\documentclass[11pt]{article}

\usepackage{amsfonts,amssymb}
\usepackage{graphicx}

\newcommand{\C}{\mathbb{C}}

\newcommand{\Hom}{\mbox{Hom}}
\newcommand{\id}{\mbox{id}}

\newtheorem{theorem}{Theorem}[section]

\newtheorem{lemma}[theorem]{Lemma}
\newtheorem{proposition}[theorem]{Proposition}

\newtheorem{definition}[theorem]{Definition}

\newtheorem{example}[theorem]{Example}

\newtheorem{remark}[theorem]{Remark}

\input{epsf.sty}

\def\DB{ {\rm DB }}
\def\K{{\mathbb K}}

\def\id{ 1 }
\def\ox{\otimes}

\begin{document}

\title{Cocycle Deformations of Algebraic Identities and $R$-matrices}

\author{
J. Scott Carter\footnote{Supported in part by
NSF Grant DMS  \#0603926.}
\\ University of South Alabama
\and
Alissa S. Crans \\
Loyola Marymount University
\and
Mohamed
Elhamdadi
\\ University of South Florida
\and
Masahico Saito\footnote{Supported in part by NSF Grant DMS
 \#0603876.}
\\ University of South Florida
}

\maketitle

\begin{abstract}
For an arbitrary identity $L=R$ between compositions
of maps
$L$ and $R$ on tensors of vector spaces $V$,
a general construction of a $2$-cocycle condition is given.
These $2$-cocycles correspond to those obtained in deformation
theories of algebras. The construction is applied to a canceling
pairings and copairings, with explicit examples with calculations.
Relations to the Kauffman bracket and knot invariants are discussed.
\end{abstract}

\section{Introduction}



The $2$-cocycle conditions
of Hochschild cohomology 
of algebras and bialgebras
are interpreted in deformations of algebras~\cite{GS}.
In other words, a map satisfying the associativity condition can be
deformed to obtain a new associative map in a larger vector space 
using $2$-cocycles.

Our motivation for this paper and other recent work \cite{CCES1,CCES2,CCEKS} comes from the fact that quandle cocycles \cite{CES, CJKLS} can be regarded as giving cocycle deformations of $R$-matrices (solutions to the Yang-Baxter equation (YBE)).  Thus it was natural to ask if this principle could be applied
to other algebraic constructions of $R$-matrices,
to construct new $R$-matrices from old via $2$-cocycle deformations. We 
have had some success in constructing
new $R$-matrices from old
using $2$-cocycle deformations.  Specifically, in \cite{CCES1}, self-distributivity was revisited
from
the
point of view of coalgebra categories,
thereby unifying Lie algebras and quandles in these categories.
Cohomology theories of Lie algebras and quandles
were given a unified definition, and deformations of $R$-matrices
were constructed.
In  \cite{CCES2}, the adjoint map of Hopf algebras, which
 corresponds to the group conjugation map, 
was 
studied from the same viewpoint. 
A cohomology theory was constructed based on  equalities satisfied
by the adjoint map that are sufficient for it to satisfy the YBE.  
Finally, in \cite{CCEKS} we presented an analog for Frobenius algebras using multiplication and comultiplication.

In the first half of this paper, we will describe a general principle of
constructing deformation $2$-cocycles
from algebraic identities (such as associativity)
that relate two (apparently) distinct tensor operators.
Then in the second half, the principle is applied to
the Kauffman bracket pairings
to construct
$2$-cocycle deformations of bracket $R$-matrices.
This is  the same approach
that we took in
\cite{CCES1,CCES2,CCEKS}. In this way, we are extending these methods to another algebraic structure.
The cup-cap pairings are among the most simple of which we can conceive. Yet that there is a deformation theory even here strikes us as interesting. In our final section,
we will discuss the knot invariants that result from these
deformed $R$-matrices.

\section{ A general construction  of  $2$-differentials}

In the deformation theory of algebras and coalgebras \cite{GS,MS},
deformation cocycles
arise as the primary obstructions to
extending
 a formal deformation of the algebraic structure.  We observe in this section that the
 cocycle conditions,
 in general,
  can be loosely described as
 ``infiltrating an algebraic condition" with an arbitrary function.
 A standard example will help illustrate the idea.
Given a $2$-cocycle
 $\phi$, the $2$-cocycle condition of an associative algebra $A$
is written as
$$ \phi (ab \otimes c)+\phi (a \otimes b)c=  \phi (a \otimes bc)+ a \phi (b \otimes c) , $$
for $a,b,c \in A$ 
or
in the multiplicative notation of linear maps
$\mu : A \otimes A \rightarrow A$
\begin{eqnarray}\label{2cocyeqn}
\ \phi (\mu \otimes \id )  + \mu ( \phi \otimes \id )
=  \phi  (\id \otimes   \mu ) + \mu (\id \otimes   \phi ) .
 \end{eqnarray}
To derive this formula,
take the
associative law
$((ab)c)=(a(bc))$, write it as
$\mu(\mu \otimes \id) = \mu ( \id \otimes \mu) $,
and put  distinct subscripts on 
  the multiplication maps
for both sides of the identity to obtain
$\mu_1 (\mu_2  \otimes \id)
 = \mu_1 ( \id \otimes \mu_2) $.
 Then take a formal sum
of
each side by replacing 
each of the maps in turn 
 by a map  $\phi$ to obtain
$$ \phi  (\mu_2  \otimes \id) + \mu_1 (\phi  \otimes \id)
 = \phi ( \id \otimes \mu_2) + \mu_1 ( \id \otimes \phi) . $$
By removing the subscripts, we obtain the $2$-cocycle condition (\ref{2cocyeqn}).
Thus we see this scheme as a cocycle $\phi$ infiltrating a formal sum of the identity.

This scenario 
has been 
generalized to a large variety of cases 
including our work in \cite{CCES1,CCES2,CCEKS}.
Our purpose in defining such generalizations is
so that we can
develop topological invariants of knots, manifolds, and knotted surfaces from cocycle conditions. The axiomatizations that we have developed are given via diagrammatic  formulations. The diagrammatic versions often lead directly to topological interpretations via equivalences such as the Pachner moves, the Yang-Baxter condition, or the tetrahedral condition for knotted surfaces.
On the other hand, the generality in which infiltration gives  chain complexes is much
 more broad than the applications that we have found. Here, we will describe the situation in a broad setting, discuss the situations for which we have found interesting results, and point to some future generalities.

\subsection{Single-term identities}

For a linear map $F: V^{\otimes p} \rightarrow V^{\otimes q}$,
we call $(|_a \otimes F \otimes |_b):
V^{\otimes a+b+p} \rightarrow V^{\otimes a+b+q}$
a map {\it expanded from $F$ by  identities}.
Here $a,b,p$, and $q$ are non-negative integers, and
the symbol
$|_x$ denotes the identity map on $x$ tensor factors of the underlying vector space $V$.
Let $ {\rm Perm}: V^{\otimes p}  \rightarrow V^{\otimes p}$
be  a composition of maps expanded from the transposition
by identities (i.e., a permutation of tensor factors, possibly
the identity), called simply a {\it permutation} on $V^{\otimes p}$.
A map written as $ (|_a \otimes F \otimes |_b)\circ  {\rm Perm}:
V^{\otimes a+b+p} \rightarrow V^{\otimes a+b+q}$
is called
{\it a map expanded from $F$ by identities and transpositions},
where Perm is a permutation on $V^{\otimes a+b+p}$.

\begin{figure}[htb]
\begin{center}
\includegraphics[width=3.5in]{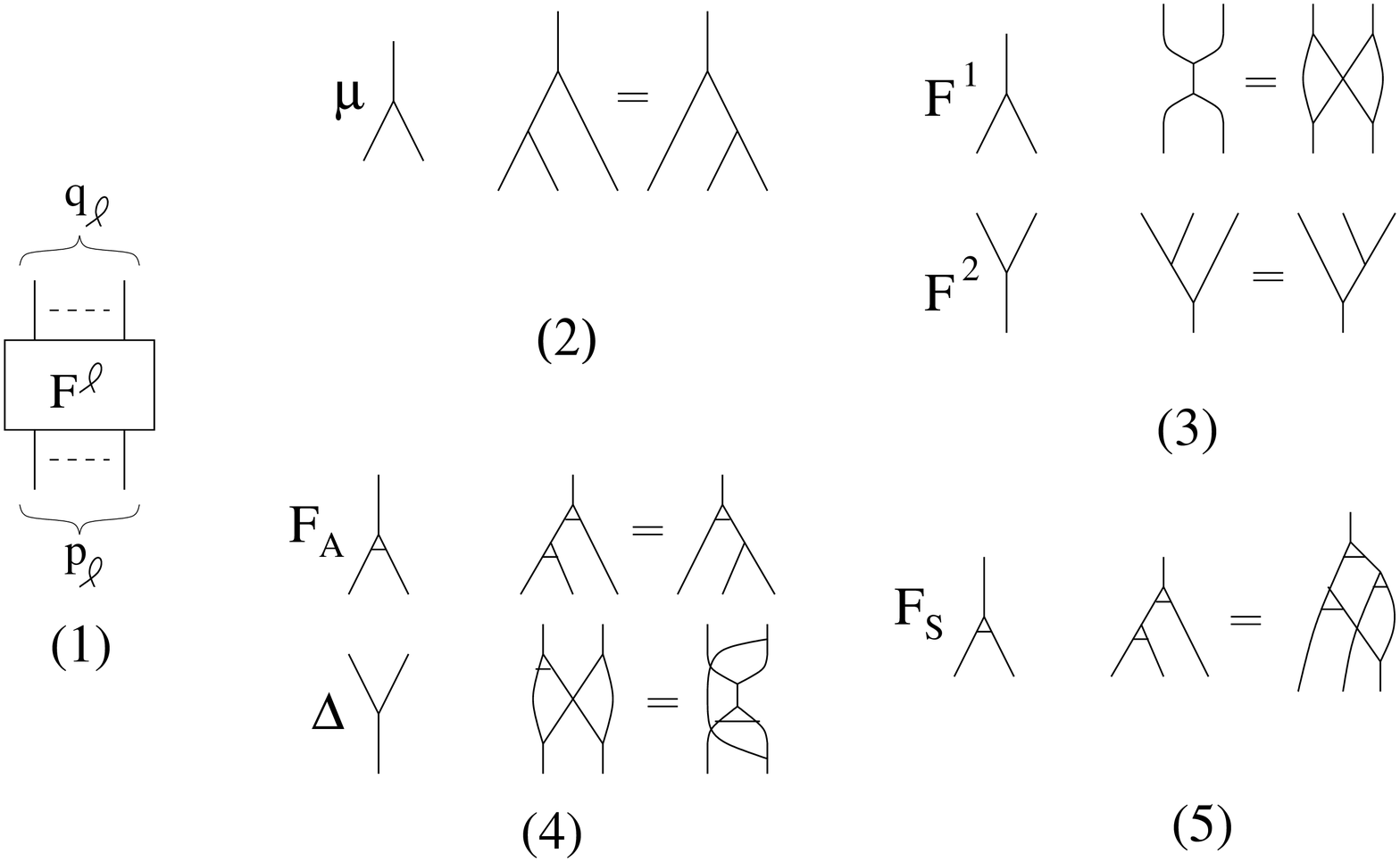}
\end{center}
\caption{Single-term identities }
\label{qpon}
\end{figure}

For a finite set of linear maps ${\cal F}=\{ F^\ell \ : \ \ell=1, \ldots, k\}$,
let $L=R$ be an
equation of linear maps
from $V^{\otimes p}$ to $V^{\otimes q}$ such that
both $L$ and $R$ are
 composites of maps expanded from $F^\ell\/$s
 ($F^\ell \in  {\cal F}$) and by identities and transpositions. We call this equation a {\it single-term identity}.
In Fig.~\ref{qpon} (1), our diagrammatic convention is depicted.
Each vertical string represents a tensor factor of $V$,
and the diagram is read from bottom to top,
so that the composition $fg$ of maps is 
represented by
the diagram of $f$ on top of that of $g$.
In (2), a multiplication
map $\mu: V \otimes V \rightarrow V$ is depicted
on the left, and the single-term identity for
associativity
is depicted
on the right.
Similarly, (2) through (5)
 depict the corresponding maps and identities
for the examples that follow.

\begin{example}\label{alg}{\rm  Suppose that $V$ is a vector space. Let $F^1: V\ox V \rightarrow V$ and let $F^2:V \rightarrow V \ox V$.
The {\it (single-term) bialgebra identities} are
\begin{eqnarray*}
F^1(F^1 \ox |) &=& F^1(| \ox F^1) , \\
(| \ox F^2)(F^2) &=& (F^2 \ox |)F^2, \\
F^2 F^1 &=&  (F^1 \ox F^1)(| \ox X \ox |)(F^2 \ox F^2)  \end{eqnarray*}
where $|$ denotes the identity
and $X$ denotes
a 
transposition
that acts, here, on  
the middle two tensorands.
These correspond to associativity, coassociativity, and compatibility, which are illustrated in (3) above.
}\end{example}

\begin{example}\label{adjoint}{\rm
 Let $V$
  be a Hopf algebra with multiplication 
  $\mu:V \ox V \rightarrow V$, comultiplication
 $\Delta:V \rightarrow V \ox V$, and antipode $S:V \rightarrow V$. Let $F_A:V\ox V \rightarrow V$ denote
 the adjoint map.
  Then $F_A= \mu(|\ox \mu)(S\ox |_2)(X \ox |)(| \ox \Delta)$  where, as before, $X$ denotes a transposition,
 and $|$ denotes the identity map. The  {\it (single-term) adjoint identities}
 consist of
the identities
\begin{eqnarray*}
F_A(F_A \ox |) &=& F_A(|\ox \mu), \\
(F_A  \ox \mu)(|\ox X \ox |)(\Delta \ox \Delta) &=& (|\ox \mu)(X\ox
|)(|\ox \Delta)(|\ox
F_A)(X\ox |)(|\ox \Delta),
\end{eqnarray*}
which are illustrated in (4) above.
}\end{example}

\begin{example}\label{shelf}{\rm
Suppose a vector space $V$ has a
cocommutative
comultiplication $\Delta:V \rightarrow V \ox V$. Then $V$ has
{\it (single-term) categorical self-distributivity}
if there is a map $F_s:V\ox V \rightarrow V$ that satisfies:
$$F_s(F_s \otimes |) = F_s (F_s \ox F_s)(|\ox X \ox |)(|\ox|\ox \Delta),$$ which is illustrated in (5) above.
In \cite{CCES1}, we also assumed that $\Delta$ satisfied coassociativity (see Example~\ref{alg}).
}\end{example}

\begin{figure}[htb]
\begin{center}
\includegraphics[width=3.5in]{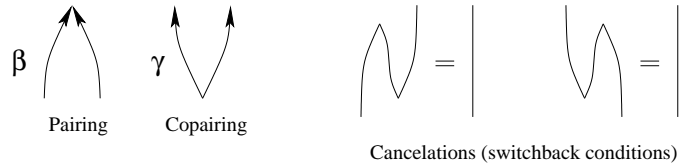}
\end{center}
\caption{Diagrams for a pairing, a copairing and their identities }
\label{pairdiag}
\end{figure}

\begin{example}\label{sw}{\rm Let $V$ denote a vector space,
$\beta:V\ox V \rightarrow \K$ denote a pairing, and $\gamma:\K\rightarrow V\ox V$ denote a copairing.
Then the
{\it (single-term) switchback identities} are
$$(| \ox\beta)(\gamma \ox |) = |, \quad {\rm and} \quad
(\beta \ox |)(| \ox \gamma) = |, $$ which are illustrated in Fig. ~\ref{pairdiag}.
}\end{example}

\subsection{Elaborate plans -- variable-distinctions of  single-term equalities}

Let $L=R$ be a single-term identity.
Let $m(\ell) $ be the number of copies of $F^{\ell}$
that appear in $L$, and write
${\cal F}_L=\{ F^\ell_j \ | \ \ell=1, \ldots, k;  \ j=1, \ldots, m(\ell ) \}$,
which
we call
{\it the distinguished variable set of $L$}.
Let ${\cal L}({\cal F}_L)$ be the
formal expression  obtained from $L$ by replacing
all $F^\ell\/$s  by $F^\ell_j\/$s with distinct $j$.
For example, for the associativity
axiom of an algebra, $\mu(\mu \otimes |)=\mu (| \otimes \mu)$,
${\cal F}=\{ \mu \}$,
$L=\mu(\mu \otimes |)$, ${\cal F}_L=\{ \mu_1, \mu_2 \}$, and
${\cal L}(\mu_1, \mu_2)= \mu_1 (\mu_2 \otimes |)$.
The same notation applies to the RHS, and
again considering the associativity axiom we have
$R=\mu(| \otimes \mu)$, and
${\cal R}(\mu_1, \mu_2)= \mu_1 (|  \otimes \mu_2)$.

\begin{definition}{\rm
The formal equality ${\cal L}({\cal F}_L)={\cal R}({\cal F}_R)$
is called the
{\it elaborate plan of the single-term identity} $L=R$.
} \end{definition}

\subsection{Infiltrations  of elaborate plans -- $2$-cocycle conditions}

To simplify the notation for substitution
of 
 a map $f=f(x_1, \ldots, x_n)$ with variables $x_i$, $i=1, \ldots, n$,
we use the notation $f(x_i=g)$ to indicate that only a single variable
$x_i$ is substituted by $g$.
That is, $f(x_i=g)=f(x_1, \ldots, x_{i-1}, g, x_{i+1}, \ldots,  x_n)$.
If multiple variables are substituted,
the set of substitution rules are indicated, as in
$f(x_i=g, x_j=h)=f(x_1, \ldots, g, \ldots, h, \ldots,  x_n)$,
where $i<j$ is assumed in this case.

\begin{definition}{\rm
Let $L=R$ be a single-term identity
among
linear
maps ${\cal F}=\{F^1,\ldots, F^k\}$.
Let ${\cal L}({\cal F}_L)={\cal R}({\cal F}_R)$ be an elaborate plan of $L=R$,
where ${\cal F}_L=\{ F^\ell_i \ | \ \ell=1, \ldots k; \  i=1, \ldots, m(\ell ) \}$
and  ${\cal F}_R=\{ F^\ell_j \ | \  \ell=1, \ldots, k;  \ j=1, \ldots, n(\ell ) \}$.
We also write it as
${\cal L}( \{ F^\ell_i \ | \ \ell=1, \ldots, k; \  i=1, \ldots, m(\ell ) \} )
= {\cal R}( \{ F^\ell_j \ | \  \ell=1, \ldots, k; \  j=1, \ldots, n(\ell ) \} ) $.
Let
$\{\phi^\ell:V^{\ox p_\ell} \rightarrow V^{\ox q_\ell}: \ell = 1, \ldots , k\}$ denote
a collection of  linear maps.
An {\it infiltration of the elaborate plan
${\cal L}({\cal F}_L)={\cal R}({\cal F}_R)$}  is the formal sum
$$ \sum_{\ell, i} {\cal L}(F^\ell_i=\phi^\ell,
F^\ell_{i'}=F^\ell (i'\neq i) )
=\sum_{\ell, j} {\cal R}(F^\ell_j=\phi^\ell,
F^\ell_{j'}=F^\ell  (j'\neq j)) .$$
Here the substitution is made as follows.
For a fixed $\ell$, there are $m(\ell)$ copies $F^\ell_i$ of $F^\ell$
in the LHS $L$. First $F^\ell_1$ is replaced by $\phi^\ell$, and all the other $F^\ell_{i'}$, $i'\neq 1$,
are replaced by the original variable $F^\ell$.
Then the second term is formally added after replacing $F^\ell_2$
by $\phi^\ell$ and other $F^\ell_{i'}$, $i'\neq 2$, are replaced by $F^\ell$.
This is repeated for all $\ell$.
} \end{definition}

\begin{definition}{\rm
We define the {\it $2$-differential} by LHS$-$RHS of
the infiltration of the elaborate plan:
$$ d^{2}
( \phi^\ell : \ell = 1, \ldots , k) =
\sum_{\ell, i} {\cal L}(F^\ell_i=\phi^\ell,
F^\ell_{i'}=F^\ell (i'\neq i) )
-\sum_{\ell, j} {\cal R}(F^\ell_j=\phi^\ell,
F^\ell_{j'}=F^\ell  (j'\neq j)
) .$$

If a set of more than one single term equalities is given,
then we define a $2$-differential for each equality.
Thus if equalities $\{ L_r= R_r \ | \ r=1 , \ldots, s \}$ are given,
denote their elaborate plans by
$\{ {\cal L}_r ({\cal F}_{L_r}) = {\cal R}_r ({\cal F}_{R_r})\ | \ r=1 , \ldots, s \}$, and
we define
$$ d^{2, r}( \phi^\ell : \ell  = 1, \ldots , k) =
\sum_{\ell, i} {\cal L}_r (F^\ell_i=\phi^\ell,
F^\ell_{i'}=F^\ell (i'\neq i) )
-\sum_{\ell, j} {\cal R}_r (F^\ell_j=\phi^\ell,
F^\ell_{j'}=F^\ell  (j'\neq j)
) .$$
In this notation, the letter $r$ specifies the equality and  $\ell$ specifies
a map. Subscripts of $F$ represent distinguished copies of a map.
}\end{definition}

\begin{example}{\rm
We infiltrate the three bialgebra identities
of Example \ref{alg} by $\phi^1$ and $\phi^2$ to obtain
$$d^{2,1}(\phi^1,\phi^2)=
\phi^1(F^1 \ox |)+ F^1(\phi^1 \ox |) -\phi^1(| \ox F^1)-F^1(| \ox \phi^1)$$
$$d^{2,2}(\phi^1,\phi^2) =
(| \ox \phi^2)(F^2)+(| \ox F^2)(\phi^2)-(\phi^2 \ox |)F^2-(F^2 \ox
|)\phi^2$$ and
$$d^{2,3}(\phi^1,\phi^2)= \phi^2 F^1 +F^2 \phi^1- (\phi^1 \ox F^1)(| \ox X \ox |)(F^2 \ox F^2)$$ $$-(F^1 \ox \phi^1)(| \ox X \ox |)(F^2 \ox F^2)-(F^1 \ox F^1)(| \ox X \ox |)(\phi^2 \ox F^2)-(F^1 \ox F^1)(| \ox X \ox |)(F^2 \ox \phi^2).$$}
\end{example}

\subsection{The first differentials and $d^2 d^1 =0$}

Let $V$ denote a finite dimensional vector space over a field $\K$, and as before,
for each $\ell=1,2, \ldots, k$, let
$F^\ell:V^{\otimes p_\ell} \rightarrow V^{\otimes q_\ell}$
denote a linear
map.
Let $f: V\rightarrow V$.

\begin{definition} {\rm
For each $\ell$, the
{\it $1$-differential} is a map
$$d^{1,\ell}: \Hom(V,V)\rightarrow \Hom(V^{\ox {p_\ell}},V^{\ox {q_{\ell}}})$$
defined by
$$ d^{1,\ell} (f) =
\sum_{i=1}^{q_\ell}
( |_{i-1} \ox f \ox |_{q_\ell - i } )  F^\ell
-
\sum_{j=1}^{p_\ell} F^\ell  ( |_{j-1} \ox f \ox |_{p_\ell - j } ) .$$
} \end{definition}

\begin{figure}[htb]
\begin{center}
\includegraphics[width=3in]{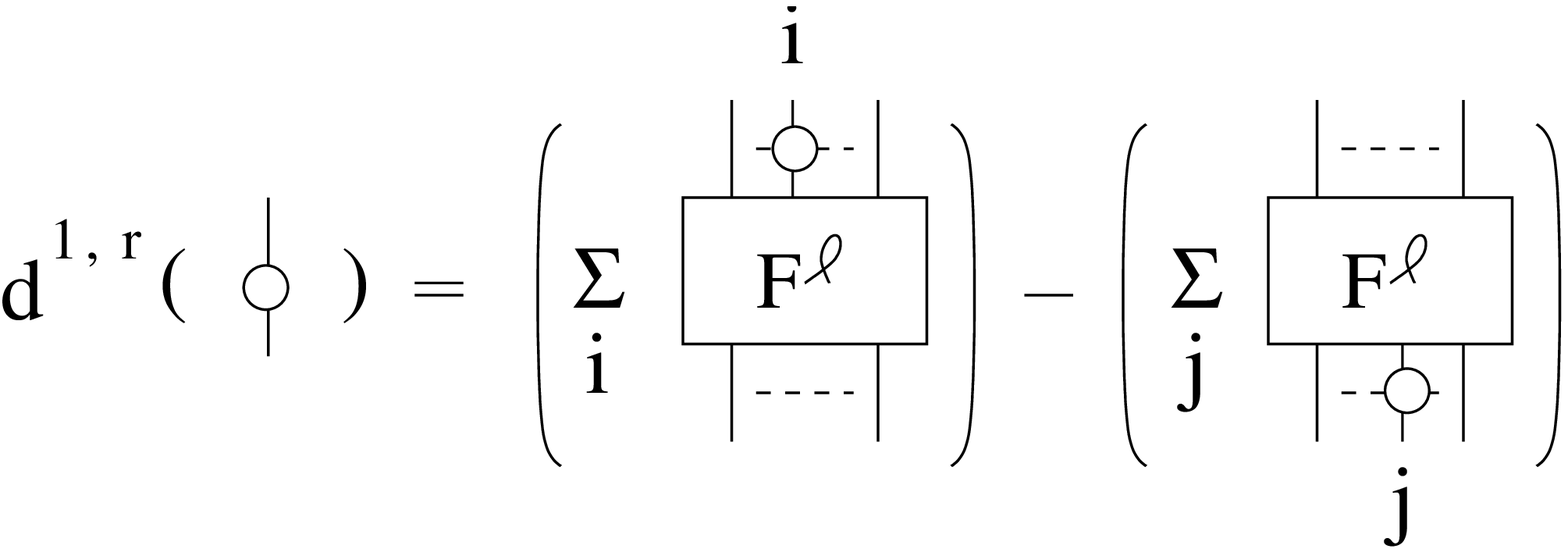}
\end{center}
\caption{The first differential $d^{1, \ell}$ }
\label{qpond1} \end{figure}

A diagrammatic representation of the $1$-differential
is depicted in Fig.~\ref{qpond1}.
The proof of the following proposition, then,
can be easily visualized by this diagram.

\begin{proposition}
Let $L=R$ be a single-term identity among
$\{ F^\ell \ | \ \ell = 1 , \ldots, k \}$,
and $d^2$ 
be the $2$-differential
of $L=R$ with variables $\{ \phi^\ell \}$.
Let $d^{1, \ell}$ be the $1$-differential for $\{ F^\ell \ | \ \ell = 1 , \ldots, k \}$.
Then
for all $\ell=1, \ldots, k$,
$$
d^2 (d^{1, 1}(f), \ldots, d^{1, k}(f)) =0 .$$
\end{proposition}
{\it Proof.\/}
Recall that each term of 
$d^2(\phi^1, \ldots, \phi^k)$
contains exactly one of $\{ \phi^1, \ldots, \phi^k \}$.
After
substituting $\phi^\ell=d^{1, \ell}(f)$,
each term
contain exactly one copy of $f$.
Hence each term, after substitution, is written as
$L_1 (|_i \ox f \ox |_j) L_2 $ or $R_1 (|_i \ox f \ox |_j) R_2 $
where $L=L_1 L_2$, $R=R_1 R_2$ are
the LHS and RHS of the given single-term identity $L=R$.
If
neither
$L_1$
nor
$L_2$ 
is 
the identity map, then
there are exactly two copies of  the term $L_1 (|_i \ox f \ox |_j) L_2 $
in 
$d^2(d^{1, 1}(f), \ldots, d^{1, k}(f))$:
one from
$\sum_{i=1}^{q_\ell}
 ( |_{i-1} \ox f \ox |_{q_\ell - i
 } )  F^\ell $, and the other from
$-
\sum_{j=1}^{p_\ell} F^\ell  ( |_{j-1} \ox f \ox |_{p_\ell - j } ) $,
and they cancel.
The same argument applies to the terms of the form
$R_1 (|_i \ox f \ox |_j) R_2 $.
After canceling, we are left with 
terms of the form
$(|_i \ox f \ox |_j) L$, $L (|_i \ox f \ox |_j) $,
$(|_i \ox f \ox |_j) R$, and $R(|_i \ox f \ox |_j) $.
The original identity implies, then, that
the terms $(|_i \ox f \ox |_j) L$ and $(|_i \ox f \ox |_j) R$
cancel, and so do the terms $L (|_i \ox f \ox |_j) $
and $R(|_i \ox f \ox |_j) $.
$\Box$

\subsection{Deformations and $2$-differentials}\label{Hochdeformsec}

The general construction defined  above is based on deformations of
algebras, and we describe the  relation for associative algebras
following \cite{MS}. Let $A$ be an associative algebra over $\K$
with multiplication $\mu: A \otimes A \rightarrow A$. A deformation
of $(A, \mu)$ is
the  
$\K[[t]]$-algebra $ (A_t, \mu_t)$, where $A_t=A
\otimes \K[[ t ]]$ with
 the multiplication  $\mu_t=\mu +t\mu_1+t^2\mu_2+ \cdots $.
  Here
  the maps $\mu_i: A \otimes A \rightarrow A$    are extended to $A_t$.
  The associativity of $\mu_t$ implies the following two equations
obtained from equating the coefficients of $t$ and
$t^2$,  respectively:
$$\mu_1(\mu(a\ox b) \otimes c)+\mu (\mu_1(a \otimes b)\ox c)
=\mu_1(a \otimes \mu (b \ox c))+\mu(a\otimes\mu_1(b \otimes c)) , $$
\vspace{-10mm}
\begin{eqnarray*}
\lefteqn{
\mu_2(
\mu(a \otimes b) 
\otimes c)+\mu_1(\mu_1(a \otimes b) \otimes c)+
\mu ( \mu_2(a \otimes b) \otimes 
c) } \\
&=&\mu_2(a \otimes \mu (b\ox c) )
+\mu_1(a \otimes \mu_1(b \otimes c))+
\mu (a\ox \mu_2(b \otimes c)) .
\end{eqnarray*}
The first equation is the Hoschschild 2-cocycle condition for $\mu_1$,
and this degree calculation explains the above general construction.

It is also known that the second equation, when   written as
\begin{eqnarray*}
\psi(a \otimes b \otimes c)
&=&\mu_1(\mu_1(a \otimes b) \otimes c)
-\mu_1(a \otimes \mu_1(b \otimes c))\\
&=&\mu (a\ox \mu_2(b \otimes c))
-\mu_2(\mu (a\ox b)  \otimes c)-
\mu ( \mu_2(a \otimes b)\ox c)
+\mu_2(a \otimes \mu (b\ox c) ),
\end{eqnarray*}
implies that $\mu_1(\mu_1(a \otimes b) \otimes c)
-\mu_1(a \otimes \mu_1(b \otimes c))$ is
a coboundary if $\mu_t$ is associative up to degree $2$.
Moreover, one can check
$\psi(a \otimes b \otimes c)
=\mu_1(\mu_1(a \otimes b) \otimes c)
-\mu_1(a \otimes \mu_1(b \otimes c))$ is a $3$-cocycle if
$\mu_1$ is a $2$-cocycle.

The argument above with respect to the deformation cocycles for an associative algebra show explicitly that these are obtained by the infiltration theory. Our prior work on cohomology of self-distributive maps, the adjoint map in a Hopf algebra, and on Frobenius algebras also can be interpreted from this infiltration theory. Pedro Lopes pointed out to us that this idea works in great generality and the discussion above is a formulation of that idea.
In the next section, we move to develop this idea in the case of 
switchback identities. 

\section{Cohomology of switchback pairs}

\subsection{Preliminaries}

It is known that a  {\it bilinear form (or pairing)}
 $\beta:  V \otimes V \rightarrow \K$
 on a vector space $V$ over a field $\K$ is
 {\it nondegenerate} if and only if
 there is a
 $\gamma : \K \rightarrow V \otimes V$
 such that
$ (\beta \otimes \id)(\id \otimes \gamma)=\id$
 and
 $ (\id \otimes \beta )(\gamma \otimes \id)=\id$.

  More generally, for a module $V$ over a unital
 ring ${\cal K}$, we
 define
  a pair $(\beta, \gamma)$, where
  $\beta \in \Hom(V^{\otimes 2}, {\cal K})$ and $\gamma \in \Hom({\cal K}, V^{\otimes 2})$,
   to be
  a {\it switchback  pair} on $V$ over ${\cal K}$ if
  they satisfy
  $ (\beta \otimes \id)(\id \otimes \gamma)=\id$
 and
 $ (\id \otimes \beta )(\gamma \otimes \id)=\id$.
 We call these conditions {\it switchback conditions}.

 Our diagrammatic conventions
 representing the
 bilinear pairing $\beta$, copairing $\gamma$,
  and the above conditions
  are depicted in Fig.~\ref{pairdiag}, from left to right, respectively.
  Parallel strings, representing tensor products
  of vector spaces, are read and oriented
  from bottom to top, when linear maps are applied.
  For a pairing, two strings merge at a single point, which
  we represent by a maximum with a corner,
  or a  cusped maximum  (instead of a smooth maximum).
  Similarly, a copairing is represented by a cornered minimum.
  Unless the orientations get confusing, it is always upward and is often abbreviated.
    These diagrammatic conventions have been used often in knot theory
  (see, for example, \cite{KP}).

\begin{example}{\bf Kauffman bracket pair.\/}\label{bracketex} {\rm
Let $V $
be a $2$-dimensional vector space over $\C$ with basis elements $x$ and $y$,
where $\beta$ is  defined on basis elements by:
$$  \beta(x \otimes x)=0, \quad \beta(x \otimes y)=i A, \quad
\beta(y \otimes x)= - i A^{-1},  \quad  \beta(y \otimes y)= 0,
 $$
where $A$ is a variable. This is the famous pairing used for the
Kauffman bracket \cite{KP}.
The corresponding copairing is defined by
$$ \gamma(1) = iA (x \otimes y) - i A^{-1} (y \otimes x) .$$

} \end{example}

\subsection{Deformations by $2$-cocycles}\label{deformsec}

Following \cite{MS} that described deformations of bialgebras, we
formulate
a deformation of
switchback pairs.
Let $(\beta, \gamma) $ be a switchback pair
on a module $V$ over a unital ring ${\cal K}$. A {\it deformation} 
of ${\cal A}=(V, \beta, \gamma) $ is
a triple 
${\cal A}_t=(V_t, \beta_t, \gamma_t)$
whose constituents are as follows: 
(1) The module $V_t=V \otimes {\cal K}[[ t ]]$
is, as indicated, the tensor product of $V$ with a formal power series. We make the identification
 $ V_t/(tV_t) \cong V$. (2) The  maps
 $(\beta_t, \gamma_t)$
of $(\beta, \gamma)$ are given by
$\beta_t= \beta + t \beta_1 + \cdots + t^n \beta_n + \cdots :
 V_t \otimes V_t \rightarrow {\cal K}$
and $\gamma_t = \gamma + t \gamma_1
+ \cdots + t^n \gamma_n + \cdots : {\cal K} \rightarrow V_t \otimes  V_t $
where $\beta_i: V \otimes V \rightarrow {\cal K} $ and
$\gamma_i :  {\cal K} \rightarrow V \otimes V  $,
$i=1, 2, \cdots$, are sequences of pairings and copairings, respectively.
Suppose $\beta$ and
$\gamma $ satisfy the
switchback conditions mod $t$,
and suppose
that  there exist $\beta_{1}: V \otimes V \rightarrow {\cal K}$ and
$\gamma_{1}: {\cal K} \rightarrow V \otimes  V $ such that
$\beta+t \beta_{1}$ and ${\gamma}+ t \gamma_{1}$
satisfy the switchback conditions mod $t^{2}$.
Define  $\xi_1, \xi_2  \in \Hom(V^{\otimes 3}, {\cal K})$ by
\begin{eqnarray*}
 ({\beta} \otimes 1)(1 \otimes {\gamma})  - \id
&=& t \xi_1
\quad {\rm mod}\  t^{2} , \label{bracoh2d1} \\
 (1 \otimes {\beta}) ({\gamma} \otimes 1) - \id
&=&  t \xi_2
 \quad {\rm mod}\  t^{2} ,  \label{bracoh2d2}
 \end{eqnarray*}
This situation describes the
{\it primary obstructions}
to formal deformations of switchback  pairs
 to be the pair of maps 
 $(\xi_1, \xi_2,)$,  as in
\cite{MS}.

For the switchback condition of  ${\beta}+t \beta_{1}$
and  ${\gamma}+t \gamma_{1}$ mod
$t^{2}$ we obtain:
\begin{eqnarray*}
 ( ({\beta} +t \beta_{1}) \otimes 1)
 (1 \otimes ({\gamma}+t \gamma_{1}))  - \id
&=& 0
\quad {\rm mod}\  t^{2} , \\
 (1 \otimes({\beta} +t\beta_{1})   )
 (   ({\gamma}+t \gamma_{1})  \otimes 1) - \id
&=&  0
 \quad {\rm mod}\  t^{2} ,
  \end{eqnarray*}
   which is equivalent by
degree calculations to:
\begin{eqnarray*}
 (d^{2,1}(\beta_{1}, \gamma_{1})= ) & &
 ( {\beta}  \otimes 1)(1 \otimes  \gamma_{1})
+ ( {\beta}_{1}  \otimes 1)(1 \otimes  \gamma)
 =\xi_1, \\
 (d^{2,2}(\beta_{1}, \gamma_{1})= ) & &
 ( \id \otimes  {\beta} )(  \gamma_{1} \otimes \id)
 +(\id  \otimes  {\beta}_{1} )( \gamma  \otimes \id)
 =\xi_2.
 \end{eqnarray*}

\begin{figure}[htb]
\begin{center}
\includegraphics[width=3in]{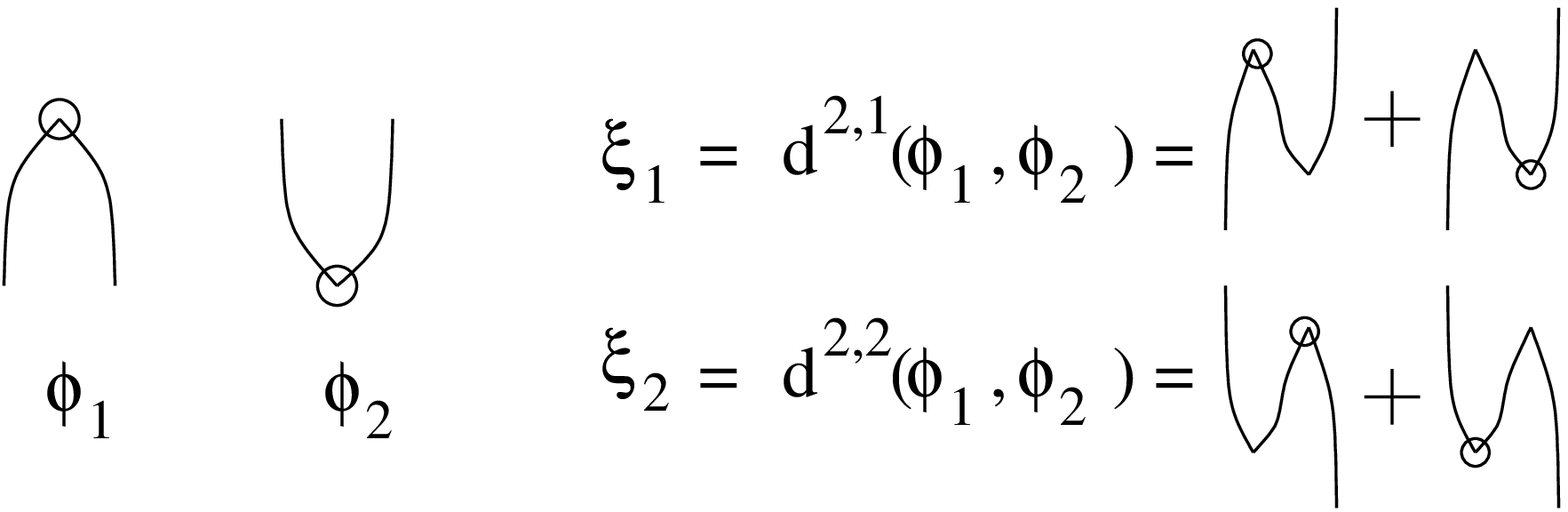}
\end{center}
\caption{The $2$-differentials}
\label{d2}
\end{figure}

Thus we make the following definition: 

\begin{definition}\label{2cocydef}
{\rm
Let $( \beta, \gamma)$
be a switchback  pair
on $V$.
Let $\phi_1 \in \Hom(V^{\otimes 2} , {\cal K})$ and
$\phi_2  \in \Hom( {\cal K}, V^{\otimes 2} )$.
Then the {\it $2$-differentials} for $(\phi_1, \phi_2)$ are
defined by
\begin{eqnarray*}
 d^{2,1}(\phi_1, \phi_2) &= &
 ( {\beta}  \otimes 1)(1 \otimes \phi_2)
+ ( \phi_1 \otimes 1)(1 \otimes  \gamma)
 , \\
 d^{2,2}(\phi_1, \phi_2) &= &
 ( \id \otimes  {\beta} )( \phi_2 \otimes \id)
 +(\id  \otimes  \phi_1 )( \gamma  \otimes \id)
 .
 \end{eqnarray*}
 If $(\phi_1, \phi_2)$ have
 vanishing $2$-differentials,
then they are called
{\it $2$-cocycles}.}
\end{definition}

In Fig.~\ref{d2}, diagrammatic representations of $(\phi_1,\phi_2)$ are
depicted
 on the left, and the $2$-differentials are depicted
  on the right.
In summary we state the following:

\begin{proposition}\label{deformprop}
{\rm (i)} The primary obstruction to deformation by $(\beta_{1}, \gamma_{1})$
of
a switchback  pair is $(\xi_1, \xi_2)$
where $\xi_1=d^{2,1}(\beta_{1}, \gamma_{1})$
and $\xi_2=d^{2,2}(\beta_{1}, \gamma_{1})$.
Hence $(\beta_{1}, \gamma_{1})$ defines a deformation
if and only if
it forms a $2$-cocycle.

\noindent
{\rm (ii)} If the primary obstruction vanishes, $(\xi_1, \xi_2)=0$,
(i.e., $(\beta_{1}, \gamma_{1})$ are $2$-cocycles),
then  the 
{\em deformation}
$(\tilde{\beta}, \tilde{\gamma})=
 ({\beta} + t  \beta_{1}, {\gamma} + t \gamma_{1})$  is a switchback  pair
 on $V_t/ (t^{2} V_t)$.
\end{proposition}

\begin{example}\label{bracketcohex}
{\rm
Let ${\cal A}=(V, \beta, \gamma)$ be
as in Example~\ref{bracketex}.
Let $\phi_1(a \otimes b)=\beta^1_{a,b} $
for basis elements $\{a,b \}=\{x,y\}$, and
$\phi_2(1)=\sum_{ \{ a,b\}=\{x,y\} } \gamma_1^{a,b} (a \otimes b)$.
Then the $2$-cocycle conditions
are formulated as:
\begin{eqnarray*}
d^{2,1}(\phi_1, \phi_2)(x)=0 & : & iA \beta_{x,x}^1 y -i A^{-1} \beta_{x,y}^1 x
+ iA \gamma_1^{y,x} x + iA \gamma_1^{y,y} y  = 0 , \\
d^{2,1}(\phi_1, \phi_2)(y)=0 & : & iA \beta_{y,x}^1 y -i A^{-1} \beta_{y,y}^1 x
- iA^{-1} \gamma_1^{x,x} x - iA^{-1} \gamma_1^{x,y} y  = 0 , \\
d^{2,2}(\phi_1, \phi_2)(x)=0 & : & - iA^{-1}\gamma_1^{x,y} x - iA^{-1}\gamma_1^{y,y} y
+ iA \beta_{y,x}^1 x -i A^{-1} \beta_{x,x}^1 y = 0 , \\
d^{2,2}(\phi_1, \phi_2)(y)=0 & : &  iA\gamma_1^{x,x} x+ iA\gamma_1^{y,x} y
+ iA \beta_{y,y}^1 x -i A^{-1} \beta_{x,y}^1 y = 0,
\end{eqnarray*}
which imply
$$\gamma_1^{y,y}=-\beta^1_{x,x}, \quad  \gamma_1^{x,x}=-\beta^1_{y,y}, \quad
  \gamma_1^{y,x}  =  A^{-2} \beta^1_{x,y}  , \quad
  \gamma_1^{x,y}  =  A^2 \beta^1_{y,x}  .
$$
Hence in total there is a $4$-dimensional solution space.
}
\end{example}

\subsection{Cohomology groups}

We discuss defining $1$-differentials
and $3$-differentials,
in relation to the above defined $2$-differentials,
and construct a chain complex in low dimensions.

\begin{definition}
{\rm
Let ${\cal A} = (V, \beta, \gamma)$ where $V$ is a module over
a unital ring ${\cal K}$, and $(\beta, \gamma)$ is a switchback pair.
Define {\it chain groups} in low dimensions as follows:
\begin{eqnarray*}
C^1( {\cal A}) &=& \Hom (V, V) , \\
C^2( {\cal A}) &=& \Hom (V^{\otimes 2}, {\cal K}) \oplus   \Hom ({\cal K}, V^{\otimes 2}) , \\
C^3( {\cal A}) &=& \Hom (V, V)_{(1)}  \oplus  \Hom (V, V)_{(2)}  , \\
C^4( {\cal A}) &=& \Hom (V^{\otimes 2}, {\cal K}) \oplus   \Hom ({\cal K}, V^{\otimes 2}) ,
\end{eqnarray*}
where the subscripts for $\Hom(V, V)$ in $C^3$ are to specify each factor.
Define differentials as follows:
\begin{eqnarray*}
d^{1,1} (\eta)  =  \beta (\eta \otimes  \id) - \beta (\id \otimes \eta) ,
& & d^{1,2} (\eta)  =  (\eta \otimes  \id)\gamma  -  (\id \otimes \eta)\gamma , \\
d^{2,1} (\phi_1, \phi_2) =
 ( {\beta}  \otimes 1)(1 \otimes \phi_2)
+ ( \phi_1 \otimes 1)(1 \otimes  \gamma)
 ,  & &
 d^{2,2}(\phi_1, \phi_2) =
 ( \id \otimes  {\beta} )( \phi_2 \otimes \id)
 +(\id  \otimes  \phi_1 )( \gamma  \otimes \id) , \\
 d^{3,1} (\xi_1, \xi_2)  =  \beta (\xi_1 \otimes  \id) - \beta (\id \otimes \xi_2) ,
& & d^{3,2} (\xi_1, \xi_2)  =  (\xi_2 \otimes  \id)\gamma  -  (\id \otimes \xi_1)\gamma ,
\end{eqnarray*}
where $\xi_1 \in \Hom (V, V)_{(1)}$
and $\xi_2 \in    \Hom (V,
V)_{(2)} $. Then further define:
\begin{eqnarray*}
 D_1(\eta) &=&  d^{1,1} (\eta) \ ( \in \Hom (V^{\otimes 2}, {\cal K}) ) \ +
 \ d^{1,2} (\eta) \ ( \in \Hom ({\cal K}, V^{\otimes 2}) ) , \\
  D_2( \phi_1, \phi_2 ) &=&  d^{2,1} (\phi_1, \phi_2) \ ( \in \Hom (V, V)_{(1)} ) \ +
 \ d^{2,2} (\phi_1, \phi_2) \ ( \in \Hom(V, V)_{(2)}  ) , \\
   D_3( \xi_1, \xi_2 ) &=&  d^{3,1} (\xi_1, \xi_2) \ ( \in \Hom (V^{\otimes 2}, {\cal K}) ) \ +
 \ d^{3,2} (\xi_1, \xi_2) \ ( \in \Hom({\cal K}, V^{\otimes 2})  ) ,
 \end{eqnarray*}
and finally,
$$ B^n ({\cal A}) = {\rm Image}(D_{n-1}), \quad
Z^n ({\cal A}) = {\rm Ker}(D_{n}), \quad
 H^n ({\cal A}) = Z^n ({\cal A}) /  B^n ({\cal A}) , $$
 for appropriate values of $n$.
 } \end{definition}

\begin{figure}[htb]
\begin{center}
\includegraphics[width=4in]{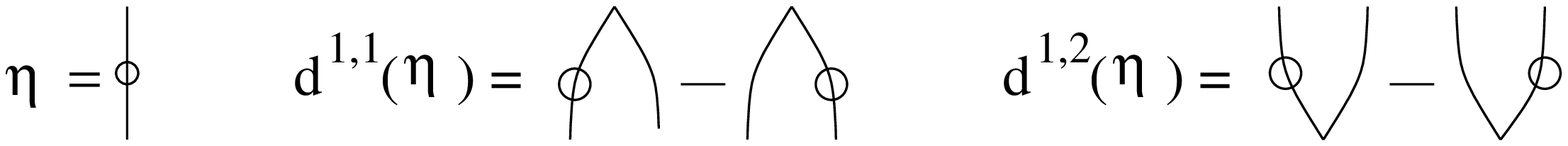}
\end{center}
\caption{The $1$-differentials}
\label{d1}
\end{figure}

Our diagrammatic conventions for representing cochains and
differentials are as follows. A
$1$-cochain $\eta \in
\Hom(V,V)$ is represented by a small white circle on a vertical
string as depicted
on the left of Fig.~\ref{d1}. The first
differentials are depicted
on the right of the figure.

\begin{figure}[htb]
\begin{center}
\includegraphics[width=2in]{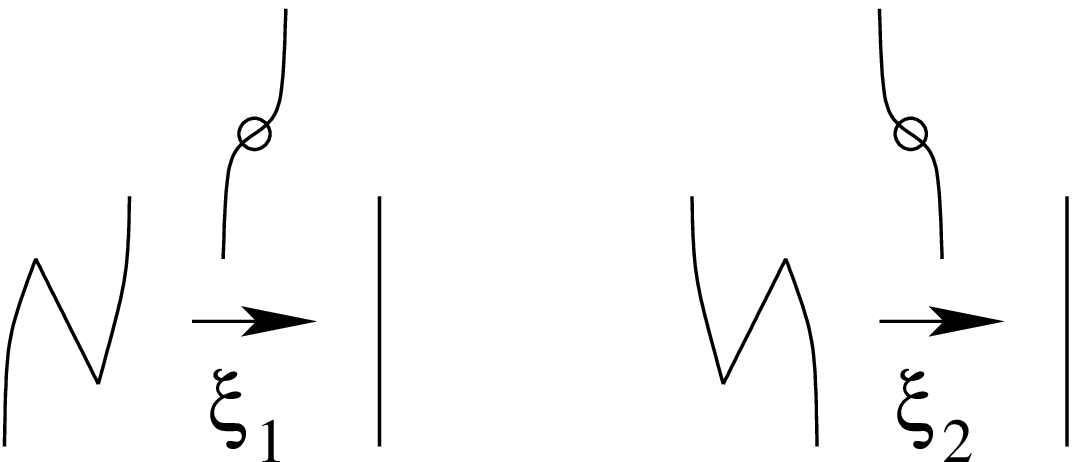}
\end{center}
\caption{Representing $3$-cocycles}
\label{cocy3}
\end{figure}

For a $3$-cochain $\xi_i \in  \Hom (V, V)_{(i)}$, $i=1,2$,
there are two aspects
 of the diagrams depicted in Fig.~\ref{cocy3}.
First, to distinguish elements in the two factors of $\Hom(V,V)$,
we use the graphs of $y=\pm x^3$ with small white circles at the origin,
respectively,
as in the figure.
When $\xi_i$, $i=1,2$,  is regarded as $d^{2,i}(\phi_1, \phi_2)$, respectively,
the graphs  $y=\pm x^3$ are regarded as cusp points
as in the figure. This is justified by the fact that
the switchback condition, when regarded as a continuous move,
corresponds to the cusp singularity of plane maps from a plane~\cite{GG}.

\begin{figure}[htb]
\begin{center}
\includegraphics[width=3in]{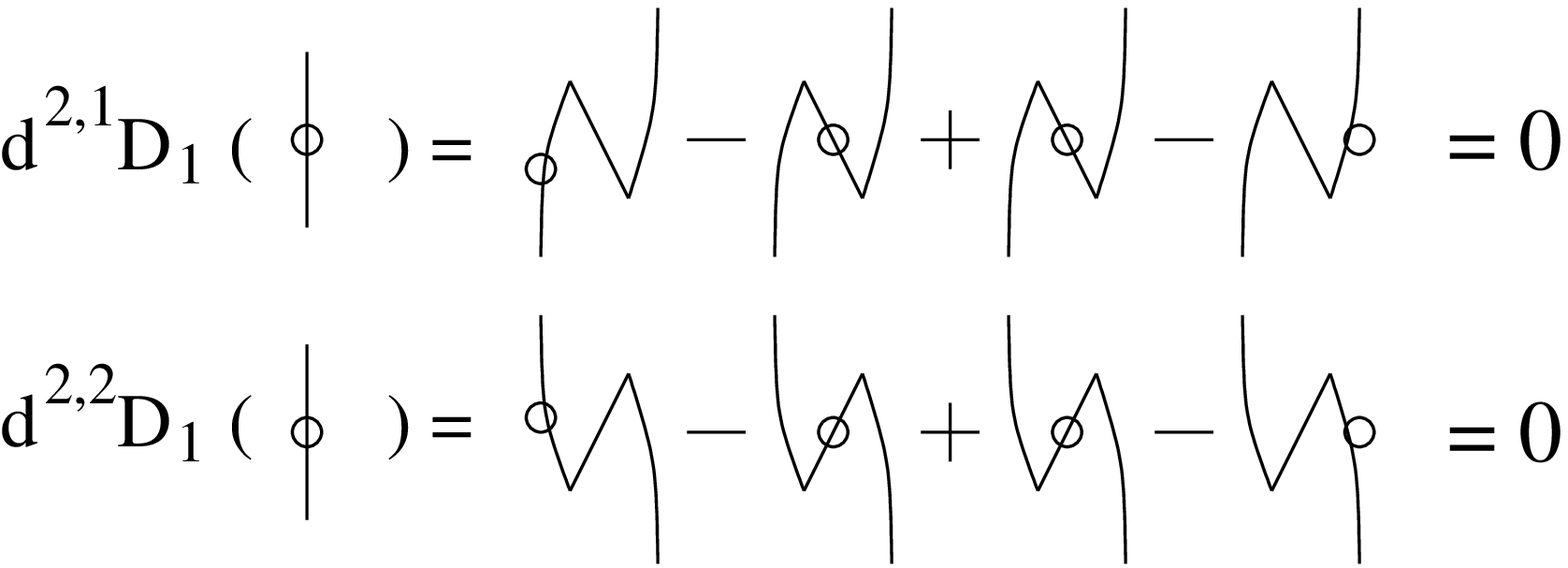}
\end{center}
\caption{$D_2 D_1=0$}
\label{d2d1}
\end{figure}

\begin{figure}[htb]
\begin{center}
\includegraphics[width=4in]{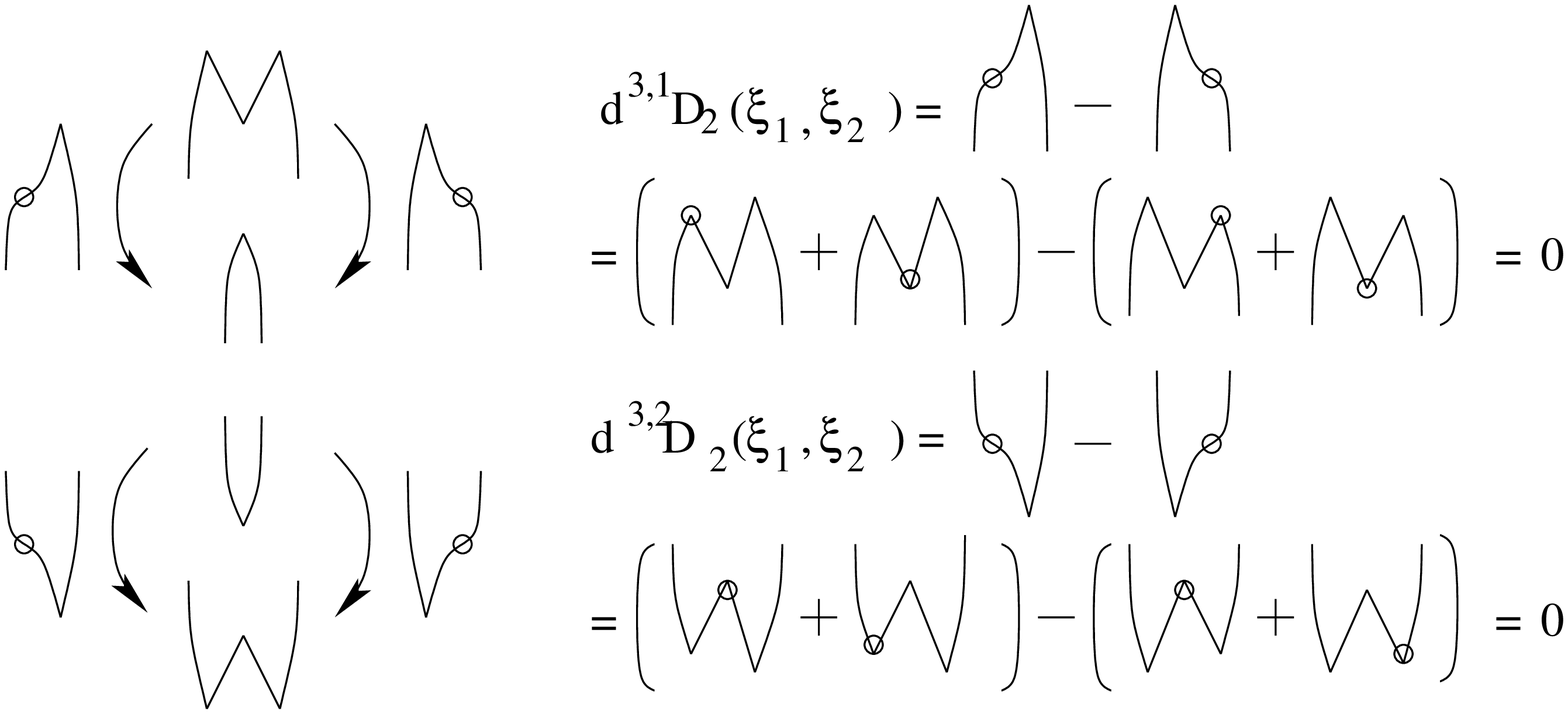}
\end{center}
\caption{$D_3 D_2=0$}
\label{d3d2}
\end{figure}

\begin{theorem}The above defined chain groups and differentials
form a chain complex:
$$ 0 \rightarrow C^1({\cal A})
\stackrel{D_1}{\rightarrow} C^2({\cal A})
\stackrel{D_2}{\rightarrow} C^3({\cal A})
\stackrel{D_3}{\rightarrow} C^4({\cal A}) . $$
\end{theorem}
{\it Proof.\/} This
 follows from direct calculations using the switchback
conditions, aided by diagrams. The fact $D_2 D_1=0$, for example, is
depicted in Fig.~\ref{d2d1}. For $D_3 D_2=0$, diagrammatic
calculations are shown in
 Fig.~\ref{d3d2}.
More specifically,
 on the left of the figure,
 we illustrate two ways to apply
switchback conditions: left first or right first, starting from the ``M''
and ``W'' shaped curves, respectively. These moves correspond to
cusps, and $\xi_i$ ($i=1,2$).
Then these diagrams are substituted by linear combinations of
other diagrams corresponding to $d^{2,i}(\phi_1, \phi_2)$, and
the terms cancel as expected.
$\Box$

\begin{example}{\rm
For ${\cal A}=(V, \beta, \gamma)$ as 
 in Examples~\ref{bracketex} and \ref{bracketcohex},
we continue and compute cohomology groups.
Let $\eta \in C^1({\cal A})$, be 
written as
$\eta(a)=\sum_{ b \in \{ x, y \} } \eta_a^b \cdot b $.
Direct calculations show that
$D_1(\eta)=0$ implies $\eta_x^x=\eta_y^y$,
 and $\eta_x^y=\eta_y^x=0$ unless $A^2+1=0$.
This implies that $Z^1({\cal A}) = H^1({\cal A}) \cong \C$ and
$B^2({\cal A}) \cong \C^3$ unless $A^2+1=0$.
Computations in Example~\ref{bracketcohex} imply that
$Z^2({\cal A})\cong \C^4$ and $B^3({\cal A})\cong \C^4$,
so that we obtain $H^2({\cal A})\cong \C$.
Let $\xi_i \in \Hom(V,V)_{(i)} \subset C^3({\cal A}) $, $i=1,2$,
be written as $\xi_i(a) = \sum_{ \{a, b \} = \{ x, y \} } {\xi_i}_a^b (b) $,
then $d^{3,1}(\xi_1, \xi_2)=0$ implies
$ {\xi_2}_x^x={\xi_1}_y^y$,
 ${\xi_2}_y^y={\xi_1}_x^x$,
${\xi_2}_x^y = -A^{-2} {\xi_1}_x^y$ and $  {\xi_2}_y^x= -A^2
{\xi_1}_y^x$.
  The second $3$-differential  $d^{3,2}(\xi_1, \xi_2)=0$ implies
the same set of
equations as the first.
Hence we obtain
 $Z^3({\cal A})\cong \C^4$ and $H^3({\cal A})=0$.

} \end{example}

\begin{remark} {\rm
The degree $2$ terms calculated 
in Section~\ref{Hochdeformsec}
for the Hochschild cohomology has the following analogue
for switchback pairs.
Let
$\beta_t= \beta + t \beta_1 + \cdots + t^n \beta_n + \cdots :
 V_t \otimes V_t \rightarrow {\cal K}$
and $\gamma_t = \gamma + t \gamma_1
+ \cdots + t^n \gamma_n + \cdots : {\cal K} \rightarrow V_t \otimes  V_t $
be formal deformations, and assume that they also satisfy the switchback
condition, which implies that the degree two terms satisfy
\begin{eqnarray}
 ( \beta  \otimes \id)(\id \otimes \gamma_2)
+ ( \beta_1  \otimes \id)(\id \otimes \gamma_1)
 +( \beta_2 \otimes \id)(\id \otimes \gamma) =0 .\label{deg2eqn}
 \end{eqnarray}
 When we set
$ \psi_1 =  ( \beta_1  \otimes \id)(\id \otimes \gamma_1) $,
and similarly
$ \psi_2 = (\id \otimes \gamma_1)  ( \beta_1  \otimes \id)$,
we obtain  two facts similar to the Hochschild case.

\noindent
(1) If the switchback relation holds up to
degree $2$, then the above Equation (\ref{deg2eqn}) holds,
and it implies that $\psi $ is a coboundary:
$$\psi_1= ( \beta_1  \otimes \id)(\id \otimes \gamma_1)
=- d^{2,1} (\beta_2, \gamma_2) = - ( \beta  \otimes \id)(\id \otimes \gamma_2)
-( \beta_2 \otimes \id)(\id \otimes \gamma). $$

\noindent
(2) The above $\psi$ is a $3$-cocycle: $D_3(\psi_1, \psi_2)=0$.
This is verified by direct calculations.

}\end{remark}

\section{Deformations of $R$-matrices by switchback pairs and Knot Invariants}

\subsection{Constructions and deformations of $R$-matrices}

In this section we present a construction of $R$-matrices from
switchback  pairs and their deformations by $2$-cocycles.

\begin{figure}[htb]
\begin{center}
\includegraphics[width=2.5in]{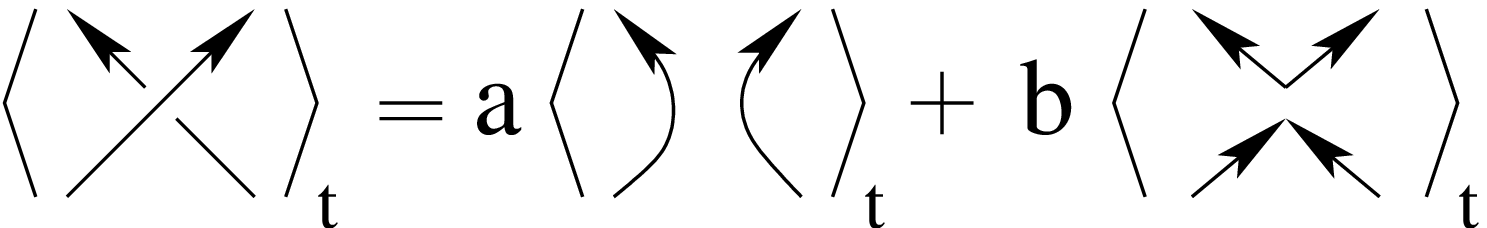}
\end{center}
\caption{The  deformed $R$-matrix}
\label{skein}
\end{figure}

\begin{lemma}\label{Rlem}
Let $(\beta, \gamma)$ be a switchback  pair on a ${\cal K}$-module $V$ as above.
 Let $\delta_0=\beta \gamma (1) \in {\cal K} $.
 Set $R=a \id + b (\gamma \beta )$ for variables
 $a$ and $b$
 taking values in ${\cal K}$.
Then $R$ is an invertible solution to the YBE  if $a$ and  $b$ are invertible
and  the equality
$
a^2 + b^2 + \delta_0 ab=0
$ holds.
\end{lemma}
{\it Proof.\/}
This lemma
 is an oriented version of the Kauffman bracket, and seems folklore.
 Direct calculations after expanding the YBE by the skein relation and
comparing the coefficients of corresponding terms gives
the above equation.
There are a few key computations of inverses,
however,
that we need to look at carefully later, and we mention these.

Set  $R'=a' \id + b' (\gamma \beta)$
for variables $a'$ and $b'$ and impose that $R'$ is the inverse of $R$ and
satisfies the YBE.
Then we obtain
additional equations,
$ a'^2 + b'^2 + \delta_0 a'b'=0$ for $R'$ to be a solution to the YBE, and
$aa'=1$, $ab'+a'b+\delta_0 b b'=0$ for $RR'=\id$.
If $a$ is invertible, we set $a'=a^{-1}$, and the equation
$a^2 + b^2 + \delta_0 ab=0$ gives
$\delta_0 b = - (a^2 + b^2) a^{-1}$, and the equation
$ab'+a'b+\delta_0 b b'=0$ becomes
$ab'+a^{-1}b  - (a^2 + b^2) a^{-1} b'=0$, i.e.,
$b(1-bb')=0$.
 With the substitutions $a'=a^{-1}$ and $b'=b^{-1}$  in
$ a'^2 + b'^2 + \delta_0 a'b'=0$, we get back $a^2 + b^2 + \delta_0 ab=0$.
Hence the result follows.
We note that the conditions are also necessary except the case $b=0$,
in which case we have  rather trivial $R$-matrices $R=a \id$ and $R^{-1}=a^{-1} \id$.
 $\Box$

The pair $(\beta, \gamma)$ in Example~\ref{bracketex},
of course,  gives the bracket,
where $a=b'=A$, $a'=b=A^{-1}$, and $\delta_0=-A^2-A^{-2}$ in the lemma.
Let $(\tilde{\beta}, \tilde{\gamma})=(\beta + t \beta_1, \gamma + t \gamma_1)$
be a deformation by $2$-cocycles $(\beta_1, \gamma_1)$.
By Proposition~\ref{deformprop},
 $(\tilde{\beta}, \tilde{\gamma})$ is a switchback  pair on $V_t/ (t^2 V_t )$
 which is regarded as a module over ${\cal K}[t]/(t^2)$.
 The above lemma applies with
 $R_t=a \id + b (\gamma + t \gamma_1)(\beta + t \beta_1)$.
 The coefficients $a,b$ and $\delta_0$ needs to be recalculated,
 and by setting $t=0$, we recover the original $R$-matrix.
 Using
 Example~\ref{bracketcohex}, we summarize this situation
 as follows for the Kauffman bracket.

 \begin{proposition}\label{Rdeformprop}
 Let $(\beta, \gamma)$ be a switchback  pair
 for the Kauffman bracket on $V$ over ${\cal K}$
  with the  $R$-matrix defined by
  $R=A \id + A^{-1} (\gamma \beta )$.
  Let  $(\beta_1, \gamma_1)$ be $2$-cocycles.
 Then the deformation
$R_t=a \id + b (\gamma + t \gamma_1)(\beta + t \beta_1)$
 is a solution to the YBE if $a^2 + b^2 + \delta_0
  ab=0$
 and $a, b$ are invertible,
 where
   \begin{eqnarray*}
\delta_0
& =& (\beta + t \beta_1)(\gamma + t \gamma_1)(1)
 =
 (-A^2 -A^{-2}) + t (iA \beta^1_{x,y} - iA^{-1} \beta^1_{y,x} ) +
 t (i A \gamma^{x,y}_1 - i A^{-1} \gamma_1^{y,x} )
 \\
 &=& (-A^2 -A^{-2})
 +  t \ [ \   i (A^2- A^{-2}) (  A^{-1} \beta^1_{x,y} + A \beta^1_{y,x} ) \   ] .
 \end{eqnarray*}
\end{proposition}

\begin{remark}{\rm
A Temperley-Lieb algebra TL$_n$ (see, for example, \cite{KP})
has generators $e_i, i=1, \ldots, n$ for a positive integer
$n$, with relations $e_i e_{i+1} e_i= e_i$,
$e_{i+1} e_{i} e_{i+1}= e_{i+1}$ for $i=1, \ldots, n-1$,
$e_i e_j=e_je_i$ for $|i-j|>1$, and $e_i^2=\delta e_i$
where $\delta \in \K$ for the coefficient field ${\K}$.
Graphically $e_i$ is represented by a pair of cup and cap as
for $\gamma \beta$.

It is well-known
[ibid] that a switchback pair $(\beta, \gamma)$
gives a representation of the  TL$_n$ by
$$e_i \mapsto \id^{\otimes (i-1)} \otimes (\gamma \beta)  \otimes  \id^{\otimes (n-i-1)}, $$
where $\delta=\delta_0=\beta\gamma(1)$.
Thus the deformation of a switchback by $2$-cocycles
$(\beta_1, \gamma_1)$
gives 
rise to
a deformation of the representation by
$$e_i \mapsto \id^{\otimes (i-1)} \otimes (\gamma + t \gamma_1)
( \beta + t \beta_1)  \otimes  \id^{\otimes (n-i-1)} $$
with $\delta=\delta_0= ( \beta + t \beta_1) (\gamma + t \gamma_1)(1)$.

}
\end{remark}

\begin{remark}{\rm
A closer inspection shows that, in fact, deformations of representations
of the Temperley-Lieb algebra TL$_n$ in the preceding remark
can be obtained from a pair $(\beta_1, \gamma_1)$
by
$$e_i \mapsto \id^{\otimes (i-1)} \otimes (\gamma + t \gamma_1)
( \beta + t \beta_1)  \otimes  \id^{\otimes (n-i-1)} $$
if they satisfy $d^{2,1}(\beta_1, \gamma_1)=d^{2,2}(\beta_1, \gamma_1)$,
which is derived from  the relations
$$e_i e_{i+1} e_i= e_i, \quad
e_{i+1} e_{i} e_{i+1}= e_{i+1} \quad \mbox{\rm  for} \quad i=1, \ldots, n-1 . $$
In the case of the 
Kauffman bracket pairings in Example~\ref{bracketex},
this condition is written as
$$\gamma_1^{x,x}=-\beta^1_{y,y}, \quad \gamma_1^{y,y}=-\beta^1_{x,x},
\quad \mbox{\rm and } \quad  \gamma_1^{x,y}- A^2 \beta^1_{y,x}
 = A^2( \gamma_1^{y,x}- A^{-2} \beta^1_{x,y} ). $$
Compare with Example~\ref{bracketcohex}.

}
\end{remark}

\subsection{Knot invariants from deformation $2$-cocycles}

For the rest of the section, we show that the cocycle deformations
of the bracket give rise to evaluations of the Jones polynomial by
truncated polynomials.

\begin{figure}[htb]
\begin{center}
\includegraphics[width=4in]{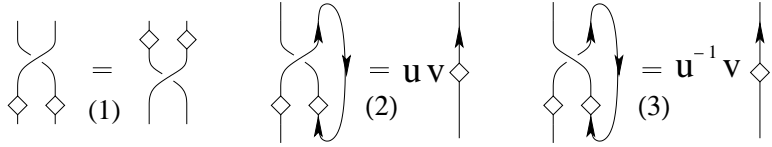}
\end{center}
\caption{Turaev's conditions }
\label{TurCond}
\end{figure}

Let $R$ be the $R$-matrix obtained by the  skein
relation in Lemma~\ref{Rlem} and its deformation
obtained in Proposition~\ref{Rdeformprop}.
We consider knot invariants obtained by  Turaev's criteria~\cite{Tur}.
For a map $f: V \otimes V \rightarrow V \otimes V$, let ${\rm Tr}_2 (f): V \rightarrow V$ denote
the map obtained from $f$ by taking the trace on the second tensor factor of
$V$. The map  ${\rm Tr}_2(f)$ is written as a composition of the
coevaluation ${\rm coev}: {\cal K} \rightarrow V \otimes V^* $ and
evaluation ${\rm ev}: V \otimes V^* \rightarrow {\cal K}$
maps by $(1 \otimes {\rm ev} ) f (1 \otimes {\rm coev} )$,
where $V^*$ denotes the dual of $V$.
These maps are defined
for basis elements $\{ v_i | i=1, \ldots , n \}$ by
${\rm coev} (1) = \sum_{i=1}^n v_i \otimes v_i^*$
and ${\rm ev}(v_i \otimes v_j)=\delta(i,j)$, where
$n$ is the dimension of $V$, and $\delta(i,j)$ is Kronecker's delta.
Diagrammatically, $f$ is represented by a box with two strings
at the top and bottom, and ${\rm Tr}_2 (f)$ is represented by the diagram of $f$ with
its 
right top and right bottom strings connected by a small loop at its right.
See the LHS of figures (2) and (3) in Fig.~\ref{TurCond}.
In this case, the right-most string representing the dual space $V^*$
is oriented downwards by convention, and the maps  ${\rm coev}$ and  ${\rm ev}$
are represented by smooth minimum and maximum, respectively,
with orientation consistently going through the maximum and minimum,
in contrast to the cusp maximum and minimum representing
pairing and copairing, with colliding orientations.

\begin{theorem}[Turaev\cite{Tur}] \label{Turthm}
Let $R: V \otimes V \rightarrow V \otimes V$ be
an (invertible)  solution to the YBE on
free module $V$ over a
commutative 
ring ${\cal K}$
with unit. 
Suppose  $\nu:  V \rightarrow V$ and $u, v \in {\cal K}$  are invertible elements that 
satisfy 
{\rm (1)} $R \circ (\nu \otimes \nu)=(\nu \otimes \nu) \circ R$,
{\rm (2)} ${\rm Tr}_2( R \circ (\nu \otimes \nu) ) = uv
\nu$, and
{\rm (3)} ${\rm Tr}_2( R^{-1} \circ (\nu \otimes \nu) ) = u^{-1}v
\nu$.
Then these maps define a link invariant via the 
closed braid $\hat{w}$ of
an $n$-braid word $w$
by
$$T_R (\hat{w})= u^{-  {\cal W}(w)} v^{- n  } {\rm Tr } ( \nu^{\otimes n}\circ  R(w) ), $$
where ${\cal W}(w)$ denotes the writhe, ${\rm Tr}$ denotes the trace,
and $R(w)$ denotes the braid group representation induced from
the $R$-matrix $R$ on $V^{\otimes n} $.
\end{theorem}

The conditions are diagrammatically depicted in Fig.~\ref{TurCond}.
The map $\nu$ is represented by a small white rhombus.

\begin{figure}[htb]
\begin{center}
\includegraphics[width=.8in]{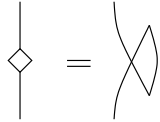}
\end{center}
\caption{The map $\nu$ }
\label{nu}
\end{figure}

\begin{figure}[htb]
\begin{center}
\includegraphics[width=4in]{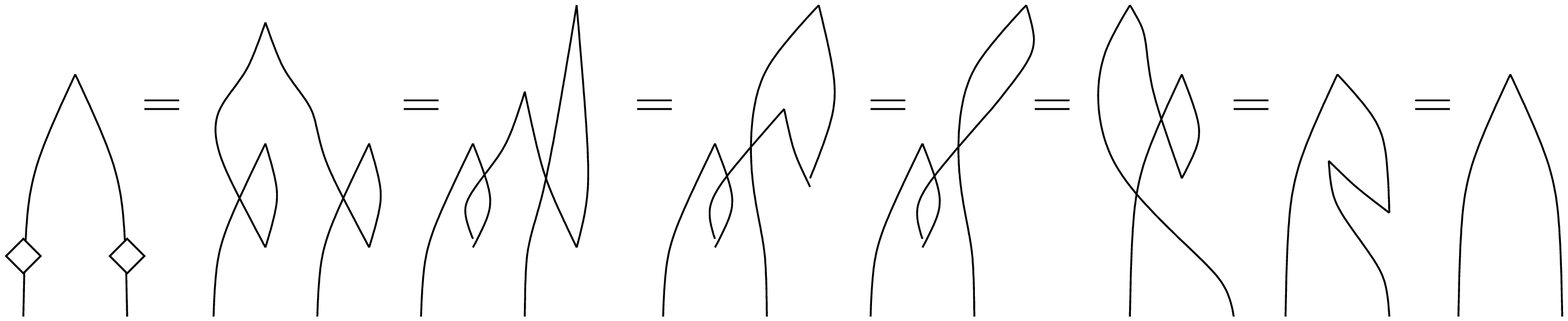}
\end{center}
\caption{Lemma~\ref{nulem} (i)}
\label{nulemma1}
\end{figure}

\begin{figure}[htb]
\begin{center}
\includegraphics[width=5in]{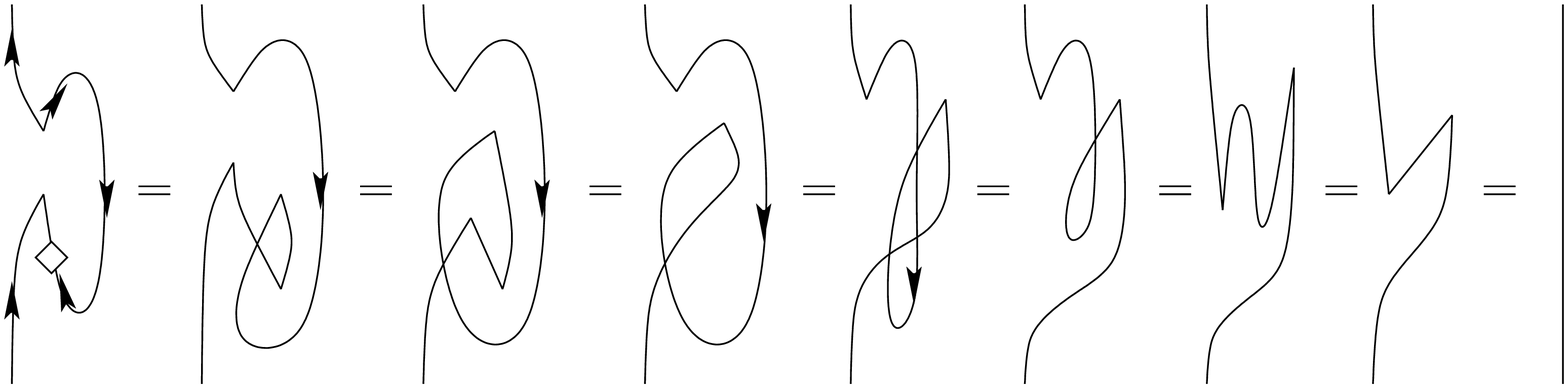}
\end{center}
\caption{Lemma~\ref{nulem} (ii)}
\label{nulemma2}
\end{figure}

To apply Turaev's construction for the $R$-matrix given in
Lemma~\ref{Rlem}, we will need the following formulas.

\begin{lemma}\label{nulem}
Let $(\beta, \gamma)$ be a switchback  pair on $V$ over ${\cal K}$.
Then for $\nu= (\id \otimes \beta ) (\tau \otimes \id)  (\id \otimes \gamma ), $
the following hold:
\quad {\rm (i) } $\beta(\nu \otimes \nu)=\beta$,  \quad
$(\nu \otimes \nu)\gamma=\gamma$. \quad
{\rm (ii) }
  ${\rm Tr}_2( \ (\gamma \otimes \id)(\beta \otimes \id)(\id \otimes \nu \otimes \id)\ )=0$.
\end{lemma}
{\it Proof.\/}
The diagram representing $\nu$ is given in Fig.~\ref{nu}.
For (i), Fig.~\ref{nulemma1} indicates a sketch of  a proof of the first equality,
and the vertical mirror images would represent a proof for the second.
For (ii), Fig.~\ref{nulemma2} shows a proof.
$\Box$

\begin{figure}[htb]
\begin{center}
\includegraphics[width=3.3in]{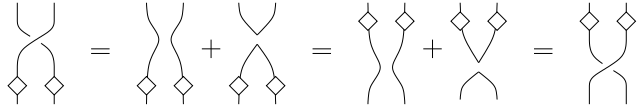}
\end{center}
\caption{Condition $(1)$ }
\label{cond1}
\end{figure}

\begin{figure}[htb]
\begin{center}
\includegraphics[width=3.3in]{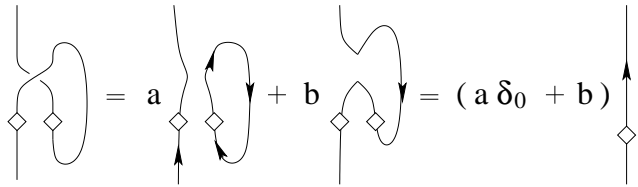}
\end{center}
\caption{Condition $(2)$ }
\label{cond2}
\end{figure}

\begin{proposition}\label{Turcondsprop}
The $R$-matrix constructed in Lemma~\ref{Rlem}
from a pair $(\beta, \gamma)$ defines a knot invariant
by Turaev's criteria with
$\nu= (\id \otimes \beta ) (\tau \otimes \id)  (\id \otimes \gamma). $

Similarly, the deformed
$R$-matrix
$R_t=a 1 + b (\gamma+t \gamma_1)(\beta + t \beta_1)$
constructed in Proposition~\ref{Rdeformprop}
defines a knot invariant
by Turaev's criteria with
$$\nu_t= (\id \otimes (\beta+ t \beta_1) ) (\tau \otimes \id)
 (\id \otimes (\gamma+t \gamma_1) ). $$
\end{proposition}
{\it Proof.\/}
This follows from Theorem~\ref{Turthm} by checking the three
conditions.
The first and the second are outlined
in Figs.~\ref{cond1}  
and  
\ref{cond2},
respectively. The third is similar to the second.
$\Box$


Let $\DB(K)$ denote the knot invariant defined
by $R_t= a \id + b ( \tilde{\gamma} \tilde{\beta} ) $ and
$R_t^{-1} = a' \id + b' ( \tilde{\gamma} \tilde{\beta} ) $
using  $\tilde{\beta}=\beta + t \beta_1$ and $\tilde{\gamma} =\gamma + t \gamma_1$
as in Proposition~\ref{Turcondsprop},
where $( \beta_1,  \gamma_1)$ are $2$-cocycles.
For the rest of the section, we compute the invariant $DB(K)$.
{}From the equalities in Fig.~\ref{TurCond},
the formulas
$u v = \delta_0 a + b $ and $u^{-1} v =  \delta_0 a' + b' $
hold,
where $aa'=1$, $bb'=1$ and
 $a^2 + b^2 + \delta_0 ab=0$. 
See the proof of Lemma~\ref{Rlem}.
The value of $\delta_0$ is given in
Proposition~\ref{Rdeformprop} as $\delta_0$.
Multiplying the skein relations
$R_t = a \id + b ( \tilde{\gamma} \tilde{\beta}) $ and
$R_t^{-1} = a' \id + b' (\tilde{\gamma} \tilde{\beta} )$
 by $u^{-1}$ and $u$, respectively, we obtain
the skein relations
$  \DB (K_+) = a u^{-1} \DB( K_0)  + bu^{-1}  \DB (K_\infty )$ and
$  \DB (K_-) = a' u \DB( K_0)  + b'u  \DB (K_\infty )$,
where $\DB(K_\infty)$ denotes the trace of a map with one crossing
replaced with $\tilde{\gamma} \tilde{\beta} $.
By eliminating the term $\DB (K_\infty )$, we obtain the relation
 $$ (b'u) \DB (K_+ ) - (b u^{-1}) \DB (K_-) = (ab' - a'b) \DB( K_0 ). $$
 Let $c=a/b$, $\ell=b^{-1}u$, and $m=ab'-a'b=c-c^{-1}$.
Then the skein relation for $\DB$ is the same as that of HOMFLYPT
polynomial. We have $\delta_0=-(c + c^{-1})$, $u=\ell b$, and $a=bc$,
and the relation for $u^2$  is simplified as follows:
\begin{eqnarray*}
 (a^{-1} \delta_0 + b^{-1})u &=& (a \delta_0 + b) u^{-1}, \\
 \ell - \ell^{-1} &=& \delta_0 ( a u^{-1} - a^{-1} u ) \\
  &=& -(c + c^{-1}) ( c \ell^{-1} - c^{-1} \ell) \\
  &=& (- c^2 \ell^{-1} + c^{-2} \ell ) + (\ell - \ell^{-1}) ,\\
  \ell^2 \ = \ c^4 , & & \ell \ =\  \pm c^2 .
 \end{eqnarray*}
This gives the skein relation of the Jones polynomial up to sign.
We summarize our calculations as follows.

\begin{proposition}\label{deforminvprop}
Let $R_t$ be the
cocycle deformation of the Kauffman bracket $R$-matrix
defined in Proposition~\ref{Rdeformprop}.
Then the knot invariant  $\DB$ defined from $R_t$ as above
is an evaluation of the Jones polynomial by a truncated polynomial.
 \end{proposition}

\subsection{Conclusion}


In this paper, we indicated that the construction of cohomology theories
via deformations
in low dimensions can be applied to broad classes of maps and identities
in variety of algebraic structures. As an example, we presented such a construction
for the Kauffman pairing and copairing, and carried out computations obtaining
 non-trivial cocycles.
 Thus the principle of cocycle deformations
 of $R$-matrices provides new solutions to the YBE.
 Probably due to the elegance of the Kauffman bracket
 and
the rigidity of the Temperley-Lieb algebra (as pointed out to us by Vaughan Jones), the resulting knot invariants are evaluations of the Jones
 polynomial by truncated polynomials.  Properties of the
  coefficients of the non-constant part,
 however, may still be of interest.
 The deformed  $R$-matrices presented in this paper, as well as in
  \cite{CJKLS, CCES1,CCES2,CCEKS},
  and  a general direction for  cocycle deformation of identities
  suggested in this paper,
 indicate unifying relations between cocycle deformations
 of algebraic systems and invariants of low dimensional knots and manifolds.

The authors would like to thank the editors of this proceedings, Vaughan Jones, and Pedro Lopes for helpful conversations.

\end{document}